\documentclass[12pt]{article}
\usepackage{graphics,amssymb,amsfonts,amsmath}
\usepackage[all,ps,cmtip]{xy} 

\usepackage{mathrsfs}
\let\mathcal\mathscr

\def\isom{\simeq}
\def\Im{\mathop{\rm Im}\nolimits}

\def\codim{\mathop{\rm codim}\nolimits}

\def\Ker{\mathop{\rm Ker}\nolimits}
\def\rank{\mathop{\rm rank}\nolimits}
\def\Sym{\mathop{\rm Sym}\nolimits}
\def\Sym{\mathbf{S}}
\def\div{\mathop{\rm div}\nolimits}

\def\iff{if and only if}

\def\lra{\longrightarrow}
\def\llra{\hbox to 12mm{\rightarrowfill}}
\def\ie{i.e.}
\def\vide{\varnothing}

\def\av{abelian variety}

\def\gndeg{nondegenerate}

\def\Pic{\mathop{\rm Pic}\nolimits}

\def\moins{\mathop{\hbox{\vrule height 3pt depth -2pt
width 5pt}\,}}



\def\N{{\bf N}}
\def\P{{\bf P}}

\def\C{{\bf C}}

\def\a{{\alpha}}

\def\cE{{\mathcal E}}
\def\cF{{\mathcal F}}
\def\cG{{\mathcal G}}

\def\cO{{\mathcal O}} 

\def\cQ{{\mathcal Q}}

\def\cS{{\mathcal S}}

\def\cV{{\mathcal V}}

\newenvironment{proof}{\trivlist \item[\hskip\labelsep{\sc
Proof.}]\rm}{\hbox
to.1pt{\hss}\hfill$\square$\bigskip\endtrivlist}

\newenvironment{proofcor}{\trivlist \item[\hskip\labelsep{\sc
Proof of the Corollary.}]\rm}{\hbox
to.1pt{\hss}\hfill$\square$\bigskip\endtrivlist}

\newtheorem{theo}{Theorem}
\newtheorem{prop}[theo]{Proposition}
\newtheorem{lemm}[theo]{Lemma}
\newtheorem{coro}[theo]{Corollary}
\newtheorem{conj}[theo]{Conjecture}

\newtheorem{rema}[theo]{Remark}
\newtheorem{remas}[theo]{Remarks}

\begin{document} 

\title{VARIETIES WITH AMPLE\\ COTANGENT BUNDLE}
\author{
Olivier Debarre}
\date{\today}
\maketitle  


Projective algebraic varieties $X$ with ample cotangent bundle have many fascinating properties: the subvarieties  of $X$ are all of general type, there are finitely many nonconstant rational maps from any fixed projective variety to $X$ (\cite{no}), if defined over the complex field,  any
entire holomorphic mapping $\C\to X$ is constant (\cite{dem}, (3.1)), if defined over a number field $K$, the set of $K$-rational points of $X$ is (conjecturally) finite (\cite{mo1}). 

However, although these varieties have been studied by several
authors (\cite{des}, \cite{mi}) and are expected to be reasonably
abundant, few concrete constructions are available.

The main result of  this article, proved
in section~\ref{s1}, is that {\em  the intersection of at least
$n/2$ sufficiently ample general hypersurfaces in an abelian variety of
dimension $n$ has ample cotangent bundle.} This answers positively a question of Lazarsfeld.
As a corollary, we obtain results about cohomology groups of sheaves of
symmetric tensors on smooth subvarieties of abelian varieties. 

In section~\ref{s2}, we discuss analogous questions
  for complete intersections in the projective space.

Finally, we present in section~\ref{s4} an unpublished result of Bogomolov
which states that a general linear section of small dimension of a product of sufficiently
many smooth projective varieties with big cotangent bundle has ample
cotangent bundle. This shows in particular that the fundamental group of a
smooth projective variety with ample cotangent bundle can be any  
group arising as the fundamental group  of a smooth projective variety.

We work over the complex numbers.

Given a vector bundle $\cE$, the projective bundle $\P(\cE)$ is the  space of
$1$-dimensional {\em quotients} of the fibers of  $\cE$. It is endowed with a
line bundle
$\cO_{\P(\cE)}(1)$. We say that
$\cE$ is
 ample (resp. nef, resp. big) if  the line bundle
$\cO_{\P(\cE)}(1)$ has the same property.
Following
\cite{som}, we say more generally that given an integer $k$, the vector bundle
$\cE$ is
$k$-ample  if, for some $m>0$, the line bundle
$\cO_{\P(\cE)}(m)$ is generated by its global sections and each fiber of the
associated map
$\P(\cE)\to
\P^N$ has dimension $\le k$. Ampleness coincide with 
$0$-ampleness.

\section{Subvarieties of abelian varieties}\label{s1}

We study the positivity properties of the cotangent bundle of a
smooth subvariety of an abelian variety $A$.

 Using a translation, we identify  the tangent
space
$T_{A,x}$ at a point $x$ of $A$  with the tangent space $T_{A,0}$ at the
origin. We begin with a classical result.

 \begin{prop}\label{prop1}
 Let
$X$ be a smooth subvariety of an \av\ $A$. The following
properties are equivalent:
 \begin{itemize}
\item[\rm (i)] the cotangent bundle $\Omega_X$ is $k$-ample;
\item[\rm (ii)] for any nonzero vector $\xi$ in $T_{A,0}$, the set $\{
x\in X\mid
\xi\in T_{X,x} \}$ has dimension $\le k$.
  \end{itemize}\end{prop}

 \begin{proof} The natural
surjection $  (\Omega_A)|_X \to\Omega_X$ induces a diagram
\begin{equation}\label{f1}
\SelectTips{cm}{}\xymatrix@M=8pt@R=20pt@C=25pt{
\P(\Omega_X)\ar@/^2pc/^f[rr] \ar@<-0.5ex>@{^{(}->}[r]^-g&
\save[]+<15.8mm,-4mm>="a"\restore
\P (\Omega_A)|_X \isom \P
(\Omega_{A,0})\times X
\ar[r]_-{p_1}&\P (\Omega_{A,0})\\& 
\save[]+<15.8mm,0mm>*\txt<8pc>{$X$}+<0mm,4mm>="b"+<-3.5mm,-3.5mm>="c"\restore\\
\ar^-{p_2}"a";"b"
\ar<8mm,-3.5mm>;"c" }
\end{equation}
\vskip-15mm
\noindent
and
$${\cal O}_{{\bf P( \Omega_X)}}(1)=g^*{\cal O}_{\P 
(\Omega_A)\vert_X}(1)=f^*{\cal O}_{\P(
\Omega_{A,0})}(1) $$
  It follows that $\Omega_X$ is $k$-ample if and only if each fiber  of $f$ has
dimension $\le k$ (\cite{som}, Corollary 1.9). The proposition follows, since the
restriction of the projection $\P( \Omega_X )\to X$ to any fiber of $f$ 
 is injective.
 \end{proof}

 \begin{remas}\upshape (1) Let $d=\dim(X)$ and $n=\dim(A)$. Since $\dim (\P
(\Omega_X))=2d-1$,  the proof of the proposition shows that the cotangent bundle
of $X$ is
$(2d-n)$-ample at best. It is always 
$d$-ample, and is   $(d-1)$-ample except if $X$ has
a nonzero vector field, which happens \iff\ $X$ is stable by translation by a
nonzero abelian subvariety (generated by the vector field).

(2) Many things can prevent the cotangent bundle of a smooth subvariety $X$ of an abelian variety $A$ from being ample. Here are two examples.
\begin{itemize}
\item If $X\supset X_1+X_2$, where $X_1$ and $X_2$ are
subvarieties of $A$ of positive dimension, then, for all $x_1$ smooth on $X_1$,
one has
$T_{X_1,x_1}\subset T_{X,x_1+x_2}$ for all $x_2\in X_2$, hence the cotangent
bundle of $X$ is not
$\bigl(\dim(X_2)-1\bigr)$-ample. In the Jacobian of a
smooth curve $C$, the cotangent bundle of any smooth $W_d(C)$ is therefore exactly
$(d-1)$-ample (although its normal bundle is ample).  
\item If $A$ is (isogenous to) a product $A_1\times A_2$ and $X_{a_2}=X\cap
(A_1\times\{a_2\})$, the cotangent bundle of $X$  is at most
$(2\dim(X_{a_2})-\dim(A_1))$-ample, because of the commutative diagram
$$\begin{matrix}
\P( \Omega_X)\vert_{(X_{a_2})_{\rm smooth}} &\stackrel{f}{\lra}&\P
(\Omega_{A,0})\\
\cup&&\cup\\
\P( \Omega_{(X_{a_2})_{\rm smooth}}) & \lra &\P
(\Omega_{A_1,0})\\
\end{matrix}
$$
In particular, if $\dim(X_{a_2})>\frac12\dim(A_1)$ for some $a_2$, the 
cotangent bundle of
$X$ cannot be ample.
\end{itemize}
%
\end{remas}  

\bigskip 
We will encounter the following situation twice: assume $\cF$ and $\cG$
are vector bundles on  a projective variety $X$  fitting into an
exact sequence
\begin{equation}\label{eq2}
0\to \cF\to  V\otimes \cO_X\to \cG\to 0
\end{equation}
where $V$ is a vector space. 

\begin{lemm}\label{le3}
In the situation above, if moreover $\rank(\cF)\ge \dim(X)   $, we have
$$\cF^*{\it\ ample}\ \Rightarrow\ \cG{\it\ nef\ and\ big}$$ 
\end{lemm}

\begin{proof}As in the proof 
of Proposition~\ref{prop1}, there is a morphism
$$f:\P(\cG) \to\P(V)
$$
 that satisfies
${\cal O}_{\P(\cG)}(1)=f^*{\cal
O}_{\P(V)}(1) $, and  $\cG$ is nef and big if and only if   $f$ is generically
finite. 

Let $d$ be the dimension of $X$, let $r$ be the dimension of $V$, let $s$ be the rank of $\cG$, and let  $G $ be 
the Grassmannian of   vector subspaces of $V^*$ of  dimension   $s$, with
tautological quotient bundle
$\cQ$ of  rank $ r -s\ge d$.
The dual of the exact sequence (\ref{eq2})   induces a map
$$\gamma:X\to G 
$$
such that $\gamma^*\cQ=\cF^*$. 

Assume that  $f\bigl( \P(\cG)\bigr)$ has dimension
$<\dim(\P(\cG))=d+s-1 $. There exists a linear subspace
$W^*$ of
$V^*$ of dimension $r-d-s+1 $ such that $\P(W)$ does not meet $f\bigl(
\P(\cG)\bigr)$.  In other words, the variety
$\gamma(X)$  does not meet the special Schubert variety
$\{\Lambda\in G\mid
 \Lambda\cap W^*\ne\{0\}\}$, whose class is $c_d(\cQ)$. It
follows that
$\gamma(X)\cdot c_d(\cQ)=0$, hence
$0=c_d(\gamma^*\cQ)=c_d(\cF^*)$.  If $\cF^*$ is ample, this
contradicts
\cite{BG}, Corollary 1.2.
\end{proof}

\subsection{Nef and big cotangent bundle}
A characterization of subvarieties   with nef
and big cotangent bundle in an abelian variety follows easily from a result of
\cite{deb}.

 \begin{prop}\label{p3}
The cotangent bundle of a smooth subvariety $X$ of an \av\ is nef
and big
\iff\
$\dim(X-X)=2\dim(X)$.
\end{prop}

 \begin{proof}The  cotangent bundle of  $X$ is nef and big \iff\ the morphism
$f$ in (\ref{f1}) is generically finite onto its image $\bigcup_{x\in
X}\P(\Omega_{X,x})$, i.e., if the latter has dimension $2\dim(X)-1$. The
proposition follows from
\cite{deb}, Theorem 2.1.
\end{proof}

The condition $\dim(X-X)=2\dim(X)$ implies of  course $2\dim(X)\le \dim(A)$. The
converse holds     if
$X$ is {\em  nondegenerate} (\cite{deb}, Proposition 1.4): this means that for
any quotient abelian variety
$\pi:A\to B$, one has either $\pi(X)=B$ or $\dim(\pi(X))=\dim(X)$.
This property holds for example for any subvariety of a {\em simple}
abelian variety.\footnote{An abelian variety $A$ is simple if
the only abelian subvarieties of $A$ are $0$ and $A$. For more about
\gndeg\ subvarieties, see
\cite{livre}, Chap.\ VIII.} 
It has
also   an interpretation in terms of   positivity of the normal bundle of
$X$.

 \begin{prop}\label{p4}
The normal bundle of a smooth nondegenerate  subvariety   of an \av\ is nef and big.
\end{prop}

 \begin{proof}Let $X$ be a   smooth    subvariety   of an \av\ $A$.
  The normal bundle $N_{X/A}$   is nef and big \iff\ the map $f'$ in
the diagram
\begin{equation}\label{f'}
\SelectTips{cm}{}\xymatrix@M=8pt@R=20pt@C=25pt{
\P(N_{X/A})\ar@/^2pc/^{f'}[rr] \ar@<-0.5ex>@{^{(}->}[r] &
\save[]+<15.3mm,-4mm>="a"\restore
\P (T_A)|_X \isom \P
(T_{A,0})\times X
\ar[r]_-{p_1}&\P (T_{A,0})\\& 
\save[]+<15.3mm,0mm>*\txt<8pc>{$X$}+<0mm,4mm>="b"+<-3.5mm,-3.5mm>="c"\restore\\
\ar^-{p_2}"a";"b"
\ar<10mm,-3.5mm>;"c" }
\end{equation}
\vskip-13mm
\noindent is  generically finite onto its image (\ie, surjective). 

To each point $p$ in the image of $f'$ corresponds a hyperplane $H_p$ 
in $T_{A,0}$ such that  $T_{X,x}\subset H_p $ for all $x $  in the image
$F_p$ in
$X$ of the fiber. This implies $T_{F_p,x}\subset H_p$ for all $x$  in
$F_p$, hence the tangent space  at the origin of the
abelian variety $K_p$ generated by 
$F_p$ is 
contained in $H_p $ (\cite{livre}, Lemme VIII.1.2).  

 Since $A$ has at most countably many
abelian subvarieties,
the abelian variety
$K_p$ is independent of the very general point $p$ in the image of $f'$. Let
$\pi:A\to B$ be the corresponding quotient. The differential of $\pi\vert_X$ is
not surjective at any point of $F_p$ since its image  is contained in the
hyperplane
$T\pi(H_p)$. By generic smoothness,  
$\pi\vert_X$ is not surjective.

If $X$ is  nondegenerate, $\pi\vert_X$ is generically finite onto
its image, hence $F_p$ is finite and $f'$ is generically finite onto its
image. It follows that $N_{X/A}$ is nef and big.
 \end{proof}

 \begin{prop}
Let $X$ be a smooth  subvariety of an abelian variety
$A$. If $
\Omega_X$ is ample, $
N_{X/A}$ is nef and big. \end{prop}

 \begin{proof}This follows from Propositions \ref{p3} and
\ref{p4} and the fact that for a \gndeg\ subvariety $X$ of $A$,  we have $\dim(X-X)=\min(2\dim(X),\dim(A))$    (\cite{deb}, Proposition 1.4).
 \end{proof}

\subsection{Ample cotangent bundle}

In this subsection, we prove that the intersection of sufficiently ample
general hypersurfaces in an abelian variety $A$ has ample cotangent bundle,
provided that its dimension be at most $ \frac12\dim(A)$.

We begin by fixing some notation. If $A$ is a smooth variety, $\partial$
a   vector field on
$A$, and $L$ a line bundle on $A$, we   define, for any section $s$ of $L$
with divisor $H$, a section
$\partial s$ of $L\vert_H$ by the requirement that for any
open set $U$ of $A$ and any trivialization
$\varphi:\cO_U\stackrel{{}_{{}_{\displaystyle\sim}}}{\lra} L\vert_U$, we have
$\partial s=\varphi(\partial(\varphi^{-1}(s)))\vert_H$ in $U\cap H$. We
denote its zero locus by
$H\cap\partial H$. We have an exact sequence
$$\begin{matrix}
H^0(A,L)&\lra& H^0(H,L\vert_H)& \lra& H^1(A,\cO_A)\\
&&\partial s&\longmapsto& \partial\smile c_1(L)\end{matrix}
$$ 
where $c_1(L) $ is considered as an element of $H^1(A,\Omega_A)$ and the cup
product is the contraction
$$H^0(A,T_A)\otimes H^1(A,\Omega_A)\lra H^1(A,\cO_A)
$$

\subsubsection{The simple case}

 We begin with the case of a simple abelian variety, where we get an explicit bound
on how ample the hypersurfaces should be.

\begin{theo}\label{va}
Let $L_1,\dots,L_c$ be very ample line bundles on a simple abelian variety
$A$ of dimension
$n$.  Consider  general divisors $H_1\in|L_1^{e_1}|,\dots,H_c\in| L_c^{e_c}|$.
If 
$e_2,\dots,e_c$ are all $>n$, the cotangent bundle of 
$H_1\cap\dots\cap H_c$ is $\max(n-2c,0)$-ample.
\end{theo}  

Since any smooth subvariety of a variety with ample cotangent bundle has the same property, it is   enough to assume $e_1,\dots,e_{[n/2]-1}$ all $>n$, and $c\ge n/2$.

\begin{proof}We need to prove that the fibers of the map $f$ in
(\ref{f1}) have dimension at most $m=\max(n-2c,0)$. This means that for $H_i$
general in $|L_i^{e_i}|$ and {\em   any} nonzero constant
vector field
$\partial$ on
$A$, the dimension of the set of points $x$ in $
X=H_1\cap\dots\cap H_c$ such that
$\partial(x)\in T_{X,x}$ is at most $m$; in other words, that
$$\dim(H_1\cap \partial H_1\cap\dots\cap
H_c  \cap  \partial H_c)\le m$$ 
It is enough to treat the case $c\le n/2$. We proceed by induction on $c$,
and assume that  the variety $Y_\partial=H_1\cap
\partial H_1\cap\dots\cap H_{c-1}  \cap  \partial H_{c-1}$ has codimension $2c-2$ in $A$ for all
nonzero  
$\partial$. Let $Y_{\partial,1},\dots,Y_{\partial,m}$  be its irreducible
components.

Let $ \cV_e(Y_{\partial,i})$ be the complement in $| L_c^e |$  of the set of
divisors
$H$ such that $(Y_{\partial,i})_{\rm red}\cap H$ is integral of
codimension $1$ in $Y_{\partial,i}$. If  
$H\notin
\cV_e(Y_{\partial,i})$, I claim that $Y_{\partial,i}\cap H\cap \partial H$ has
codimension $2$ in $Y_{\partial,i}$. Indeed, let $s\in H^0(A, L_c^e)$ define
$H$ and set $Y=(Y_{\partial,i})_{\rm red}$, with normalization  $\nu:\widehat Y\to Y$. The scheme $Y \cap H\cap
\partial H$ is the zero set in
$Y \cap H$ of the section
$\partial s$ defined above. In  the commutative diagram  
\begin{equation}\label{diaa}
\SelectTips{cm}{}\xymatrix@M=8pt@R=-8pt@C=15pt{
H^0( H, L_c^e\vert_H)\ar[r]\ar@/_2pc/[ddd]_{\nu^*} 
& H^1(A,\cO_A)\quad\ar@/^2pc/[ddd]^{\rho}\\
\quad\quad\quad\quad\partial s\ar@{|->}[r]&\partial\smile
e\, c_1(L_c)\quad\quad\quad\\
\\
H^0(\nu^*H, \nu^*L_c^e )\ar[r]&H^1(\widehat Y,\cO_{\widehat Y} )\quad }
\end{equation}
the restriction $\rho$ is injective because
$Y $ generates
$A$.

It follows that for $H\notin \cV_e(Y)$, the zero locus of 
the nonzero section 
 $\nu^*(\partial s )$ in the integral subscheme
$\nu^* H$ of $\widehat Y $  has  dimension $\le \dim(Y)-2$. For
  $H\notin \cV_e(Y_\partial)= \bigcup_{i=1}^m
\cV_e(Y_{\partial,i})$, the scheme $Y_\partial\cap H\cap \partial H$ therefore has
codimension $2c$ in $A$. Thus, for
$H_c\notin \bigcup_{[\partial]\in
\P(\Omega_{A,0})}\cV_e(Y_\partial)$, the intersection
$$ H_1\cap \partial H_1\cap\dots\cap
H_c  \cap  \partial H_c $$
has codimension $2c$ in $A$ for all nonzero constant
vector field 
$\partial$ on
$A$ (note that when $c=1$, there is no condition on $H_1$). Lemma
\ref{le} below shows that  
$\cV_e(Y_\partial)$ has codimension at least $e-1$ in $|L_c^e|$. For
$e>n$, the union $\bigcup_{[\partial]\in
\P(\Omega_{A,0})}\cV_e(Y_\partial)$ is therefore not the whole of
$|L_c^e|$ and the theorem follows.
\end{proof}

The following lemma used in the proof above is an easy consequence of \cite{ben}, Th\'eor\`eme 0.5.

\begin{lemm}[O. Benoist] \label{le}
 Let  $Y$ be an integral subscheme of $\P^n$ of dimension at least $2$, let $\nu:\widehat Y\to Y$ be its normalization, and let
$\cV_{e,n}$ be the projective space of  hypersurfaces of degree $e$ in
$\P^n$. The codimension of the complement $\cV_e(Y)$ of
$$ \{ F\in \cV_{e,n}\mid  \nu^*F\ \hbox{is integral of
codimension $1$ in $\widehat Y$ }\}
$$
in $\cV_{e,n}$ is at least $e-1$. 
\end{lemm}

\subsubsection{The general case}

A variant of the same proof works for any abelian variety, but we lose control of the
explicit lower bounds on  $e_2,\dots,e_c$.
  
\begin{theo}\label{va2}
Let $L_1,\dots,L_c$ be very ample line bundles on an abelian variety
$A$ of dimension
$n$. For $e_2,\dots,e_c$ large and divisible enough positive integers and 
general divisors $H_1\in|L_1^{e_1}|,\dots,H_c\in|L_c^{e_c}|$, the cotangent
bundle of 
$H_1\cap\dots\cap H_c$ is $\max(n-2c,0)$-ample.
\end{theo}  

 Let us be more precise about the condition on the $e_i$.  What
we mean is that there exist for each $i\in\{1,\dots,c-1\}$ a function 
$\delta_i:\N^i\to\N^*$ such that the
conclusion of the theorem holds if
\begin{eqnarray}\label{*}
 e_2=e'_2\delta_1(e_1),\
e_3=e'_3\delta_2(e_1,e_2)&\!\!\!\!,\dots,e_c=
e'_c\delta_{c-1}(e_1,\dots,e_{c-1})\nonumber\\
&\hbox{with }  e'_2,\dots,e'_c>n
\end{eqnarray}

\begin{proof}
We keep the setting and notation of the proof of Theorem \ref{va}. Everything goes
through except when, in diagram
(\ref{diaa}), $\rho(\partial\smile\, c_1(L_c))=0$. In this case, let
 $A''$ be the abelian subvariety  of $A$ generated by $Y$
 and let $A'$ be its complement with respect to $L_c$, so that
the addition
$$\pi:A'\times A''\to A$$
is an isogeny and $\pi^*L_c\isom L_c\vert_{A'}\boxtimes L_c\vert_{A''} $. We have
$Y=a'+Y''$, with
$a'\in A'$, $Y''\subset A''$, and
$\partial \in H^0(A',T_{A'})$. In particular, we have an injection
$$H^0(A,L_c)\stackrel{\pi^*}{\hookrightarrow}H^0(A',L_c\vert_{A'})\otimes
H^0(A'',L_c\vert_{A''}) 
$$ It is however difficult to identify in a manner useful for our purposes
the sections of 
$L_c$ inside this tensor product. Instead, we use a trick that will
unfortunately force us to lose any control of the numbers involved.

 The trick goes as follows. The kernel of $\pi$, being finite, is
contained in the group of
$r$-torsion points of $A'\times A''$ for some positive integer $r$. Multiplication
by
$r$ factors as
$$A'\times A''\stackrel{\pi}{\lra} A\stackrel{\pi'}{\lra}A'\times A''$$
and ${\pi'}^*( L_c\vert_{A'}\boxtimes L_c\vert_{A''}) $ is some power
$L_c^{e_0}$ of
$L_c$. Sections of $ L_c^{ e_0} $ that come from
$H^0(A', L_c \vert_{A'})\otimes H^0(A'', L_c \vert_{A''}) $ 
induce a morphism from $A$ to some projective space that factors through
$\pi'$ and embeds
$A'\times A'' $.

We will consider sections of $L^{ee_0}$ of the type
${\pi'}^*s$, with $s\in  
H^0(A', L_c^e\vert_{A'})\otimes H^0(A'', L_c^e\vert_{A''})  $.
We have $\pi'(Y)=\{ra'\}\times rY''$; let $\nu:\widehat{\pi'(Y)}\to \pi'(Y)$ be the normalization.
If the divisor $H$ of $s$ on $A'\times A''$  corresponds to a
degree
$e$ hypersurface outside of  $ \cV_e(\pi'(Y))$,  the pull-back
$\nu^* H $ is integral of codimension  1  in
$\widehat{\pi'(Y)}$.

Fix a basis
$(s_1'',\dots,s''_d)$  for
$H^0(A'',L_c^e\vert_{A''})$ and write 
$s= \sum_{i=1}^d s'_i\otimes s''_i $, so that
\begin{eqnarray*}
\nu^* H &=&\div \Bigl(\nu^*  \sum_{i=1}^d s'_i(ra')s''_i  \Bigr)\\
\nu^*( H\cap \partial H)&=&  \div \Bigl(\nu^* \sum_{i=1}^d
s'_i(ra')s''_i  \Bigr)\cap \div \Bigl(\nu^* 
\sum_{i=1}^d
\partial s'_i(ra')s''_i  \Bigr)
\end{eqnarray*}
Since $\nu^* H $ is integral, $\pi'(Y)\cap H\cap \partial H$   has
codimension $2$ in $\pi'(Y)$ (hence 
$ Y \cap {\pi'}^{-1}(H)\cap \partial {\pi'}^{-1}(H)$   has
codimension $2$ in $ Y $) unless, for some complex number
$\lambda$, the section
$ \sum_{i=1}^d
(\lambda s'_i+\partial s'_i)(ra')s''_i $ of $L_c^e\vert_{A''}$ vanishes on
$rY''$. In other words, if we let
$$\Gamma_{rY''}=\{(a_1,\dots,a_d)\in\C^d\mid \sum_{i=1}^da_is''_i\
\ \hbox{vanishes on $rY''$} \}$$ 
and
$$
M_\partial=\begin{pmatrix}
s'_1(ra')&\cdots&s'_d(ra')\\
\partial s'_1(ra')&\cdots&\partial s'_d(ra')\\
\end{pmatrix}
$$
we have $(\lambda,1)\cdot M_\partial\in
\Gamma_{rY''}$. Now we may pick any 
collection $(s'_1,\dots,s'_d)$ we like. Fix one such that the corresponding matrix
$M_\partial$ has rank $2$ for all nonzero $\partial$ and apply a square matrix $N$
of size $\dim(A')$. The condition is now that the composition
$$\Im({}^tM_\partial )\subset \C^d\stackrel{{}^tN}{\lra}  \C^d\lra
\C^d/\Gamma_{rY''}
$$
is {\em not} injective, that is, 
\begin{itemize}
\item
either  ${}^tN\cdot \Im
({}^tM_\partial)\cap\Gamma_{rY''}\ne\{0\}$, which imposes
$\codim(\Gamma_{rY''})-1$ conditions on $N$;
\item
or   
$\Ker({}^tN)\cap \Im ({}^tM_\partial) \ne\{0\}$, which imposes $d-1$
conditions on $N$.
\end{itemize}

  The ``bad'' locus for $H$
corresponds to the space of matrices $N$ that satisfies either one of these properties for
some nonzero $\partial\in H^0(A',T_{A'})$. Since, on the one hand $d=h^0(A'', L_c^e\vert_{A''})>
e$ and, on the other hand, the codimension of $\Gamma_{rY''}$ is the rank
of the linear map
$  H^0(A'', L_c^e\vert_{A''})\to H^0(rY'', L_c^e\vert_{rY''})$, which is  
$ > e$, the codimension of the  ``bad'' locus is at least $e-\dim(A')+2$.

This means that for $A''$ (hence $A'$) fixed, $e>n$, and $H$
general in $|L_c^{ee_0}|$, for any component $Y$
of $Y_\partial$ that spans (as a group)
$A''$, the intersection $Y \cap H\cap \partial H$ has codimension $2$ in
$Y$ for all nonzero $\partial$ in $H^0(A,T_A)$.

Since $A$ has at most countably many abelian subvarieties, there are only
finitely many different abelian subvarieties spanned by  components of
$Y_\partial=H_1\cap
\partial H_1\cap\dots\cap H_{c-1}  \cap  \partial H_{c-1}$ for $ H_1,\dots,
H_{c-1}$ general in $|L_1^{e_1}|,\dots,|L_{c-1}^{e_{c-1}}|$   as $\partial$
runs through the nonzero elements of $H^0(A,T_A)$. Therefore, for some
positive  integer
$\delta$, any $e>n$,   and
$H$ general in $| L_c^{e\delta}|$, the intersection $Y_\partial \cap H\cap
\partial H$ has codimension $2$ in
$Y_\partial$ for all nonzero $\partial\in H^0(A,T_A)$. This proves our claim
by induction on $c$ hence the theorem.
\end{proof}

\subsubsection{The four-dimensional  case}

In case the ambiant abelian variety has dimension $4$, we can make  the
numerical conditions in Theorem \ref{va2} explicit.
  
\begin{theo}\label{th4}
Let $L_1$ and $L_2$ be  line bundles on an abelian fourfold
$A$, with $L_1$ ample and $L_2$ very ample. For  $e_1 $ and $e_2\ge 5$,
and  $H_1\in |L_1^{e_1}|$ and $H_2\in |L_2^{e_2}|$ general,  the surface
$H_1\cap H_2$ has ample cotangent bundle.
\end{theo}

\begin{proof}I claim that  for $H_1$ general in
$|L_1^{e_1}|$, the scheme $Y_\partial=H_1\cap \partial H_1$ is an integral
surface for each nonzero   vector field
$\partial$ on
$A$. Granting the claim for the moment and using the notation of the proof of
Theorem \ref{va}, the scheme $H_1\cap \partial H_1\cap H_2$ is then, for
$H_2\in |L_2^{e_2}|\moins \cV_{e_2}(Y_\partial)$, an integral curve that
generates $A$ since its class is $e_2H_1^2H_2$. The argument of the
proof of Theorem \ref{va} applies in this case to prove  that  $ H_1\cap
\partial H_1 \cap H_2  \cap  \partial H_2 $ is finite. Taking $H_2$ outside
$  \bigcup_{[\partial]\in
\P(\Omega_{A,0})}\cV_{e_2}(Y_\partial)$ (which is possible by Lemma
\ref{le} since $e_2>4$), the intersection
$$ H_1\cap \partial H_1\cap H_2 \cap  \partial H_2 $$ is finite for all
nonzero   vector fields 
$\partial$, which is what we need. The theorem therefore follows from the
claim, proved in the next lemma.\end{proof}

 \begin{lemm}\label{lee}
Let $A$ be an abelian variety of dimension at least $4$ and let $L$ be an ample
divisor on $A$. For $e\ge5$ and $H$ general in $|L^e|$, the
scheme $ H\cap\partial H$ is  
integral for all nonzero $\partial\in H^0(A,T_A)$.
\end{lemm}

\begin{proof}Assume to the contrary that for some smooth $H\in|L^e|$, we have 
$H\cap
\partial H=D'_1+D'_2$, where $D'_1$ and $D'_2$ are effective nonzero Cartier
divisors in
$H$. We follow  \cite{bd},  proposition 1.6:  since $\dim(H)\ge3$, there exist by
 the Lefschetz Theorem
divisors
$D_1$ and
$D_2$ on $A$ such that
$D_1+D_2\equiv H$ and $D_i\vert_H\equiv D'_i$. Since $D'_i$ is effective, the long exact sequence in cohomology associated with the exact sequence
$$0\to\cO_A(D_i-H)\to \cO_A(D_i)\to\cO_H(D'_i)\to 0
$$
shows that, for each $i\in\{1,2\}$, either $H^0(A,D_i)\ne0$ or  $H^1(A,D_i-H)\ne0$.
The case where both $H^1(A,D_1-H)$ and $H^1(A,D_2-H)$ are zero is impossible,
since we would then have a section of $L^e$ with divisor $
H\cap
\partial H$ on $H$. The case where both $H^1(A,D_1-H)$ and $H^1(A,D_2-H)$ are
nonzero is impossible as in  {\em loc.~cit.} because $\dim(A)\ge3$.

 So we may assume 
$H^1(A,D_2-H)\ne0$ and $H^1(A, D_1-H)=0$, and take $D_1$ effective
such that $D_1\cap H=D'_1$.

 As in  {\em loc.~cit.}, $A$ contains an elliptic curve $E$ such that, if
$B$ is the neutral component of the kernel of the composed morphism
$$A\stackrel{\phi_H}{\lra}  \Pic^0(A)\to\Pic^0(E) 
$$
the addition map $\pi:E\times B\to A$ is
an isogeny, $\partial $ is tangent to $E$, and $\pi^*(D_1)=p_1^*(D_E)  $ for
some effective divisor $D_E$ on $E$. Pick a basis
$(t_1,\dots,t_d)$ for
$H^0(B,L^e\vert_B)$ and a section $s$ of $L^e$ with divisor $H$, and  write
$$\pi^*s=\sum_{i=1}^ds_i\otimes t_i$$ with $s_1,\dots,s_d\in
H^0(E,L^e\vert_E)$, so that $\pi^{-1}\bigl( H\cap\partial H\bigr)$  is defined by
$$  \sum_{i=1}^ds_i\otimes
t_i  = \sum_{i=1}^d\partial s_i\otimes t_i =0$$
Since $D'_1=H\cap D_1$ is contained in $ H\cap\partial H$,  for every point
$x$ of the support of
$D_E$, we have
$$ \div\Bigl(\sum_{i=1}^ds_i(x)  t_i\Bigr)\subset
\div\Bigl(\sum_{i=1}^d\partial s_i(x)  t_i\Bigr)\subset B$$
Since these two divisors belong to the same linear series $|L^e\vert_B|$  on $B$,
they must be equal and
$$\rank 
\begin{pmatrix}s_1(x)&\cdots&s_d(x)\\\partial s_1(x)&\cdots&\partial s_d(x)
\end{pmatrix}\le 1
$$Since $H$ is irreducible, the sections $s_1,\dots,s_d$ have no common zero and
the morphism $\psi_H:E\to\P^{d-1}$ that they define is ramified at $x$.

The vector subspace of $H^0(E,L^e\vert_E) $ generated by $s_1,\dots,s_d $
only depends on $s$, not on the choice of the basis $(t_1,\dots,t_d) $. If
$b_1,\dots,b_d$ are general points of $B$, it is also generated by
$s({}\cdot+b_1),\dots,s({}\cdot+b_d)$ and
$$\rank 
\begin{pmatrix}s(x+b_1)&\cdots&s(x+b_d)\\ \partial
s(x+b_1)&\cdots&\partial s(x+b_d)
\end{pmatrix}\le 1
$$ Assume now that the conclusion of the lemma fails for {\em general} $H$
(and $s$). The point $x$ varies with $s$, but remains constant for $s$ in
a hypersurface $H_x$ of $ H^0(X,L^e)$. If $s$ is in
$$H_x'=H_x\cap\{ t\in H^0(X,L^e)\mid t(x+b_1)=t(x+b_2)=0\}
$$ it also satisfies  $\partial s (x+b_1)=\partial s (x+b_2)=0$. Since
$H'_x$ has codimension at most $3$ in $ H^0(X,L^e)$, this means that $L^e$
is not
$3$-jet ample and contradicts Theorem 1 of \cite{BS}
(see also \cite{PP}): the lemma is proved.\end{proof}

\begin{remas}\upshape
(1) Let $A$ be an abelian fourfold   that contains no elliptic curves. The proof of Lemma \ref{lee} shows that for {\em any} smooth ample hypersurface $H$ in $A$ and any nonzero $\partial\in H^0(A,T_A) $, the scheme $H\cap \partial H$ is integral. It follows that for $L$ very ample, $e\ge 5$, and $H'\in |L^e|$ general,  the surface
$H\cap H'$ has ample cotangent bundle (this is   a small improvement of Theorem \ref{va}).

(2) It is proved in \cite{di} that on a {\em general}  principaly polarized abelian fourfold, the intersection of two general translates of a theta divisor   is a smooth surface with ample cotangent bundle. 
\end{remas}

\subsection{Cohomology of symmetric tensors}\label{ss13}

Let   $X$ be a smooth
subvariety of an abelian variety. We are interested in the cohomology
groups of the vector bundles $\Sym^r\Omega_X$.

\begin{prop}\label{p15}
Let $A$ be an abelian variety of dimension $n$ and let $X$ be a smooth
subvariety of codimension $c$ of  $A$ with ample normal bundle. For
$r\ge0$, the restriction
$$H^q(A,\Sym^r\Omega_A)\lra H^q(X,\Sym^r\Omega_X ) 
$$ is bijective  for $q<n-2c$ and injective  for $q=n-2c$.\footnote{For
the case $q=0$, Bogomolov gave in \cite{bog4} a very nice proof  that goes
as follows. Arguing as in the proof of Proposition
\ref{p3}, we find that the morphism $f$ of (\ref{f1}) is surjective
whenever
$X-X=A$. Any fiber of  
$f$ is isomorphic to its projection to $X$, which is the zero
locus of a section of $N_{A/X}$. It follows that when $N_{A/X}$ is ample and
$c<n-c$, the fibers of $f$ are connected, hence $f_*\cO_{\P(
\Omega_X)}(r)\isom
\cO_{\P (\Omega_{A,0})}(r)$, from which we get, for all $r\ge0$,
\begin{eqnarray*} H^0(X,\Sym^r\Omega_X)&\isom& H^0(\P( \Omega_X),\cO_{\P(
\Omega_X)}(r))\\ &\isom&  H^0(\P (\Omega_{A,0}),\cO_{\P (\Omega_{A,0})}(r))
 \isom H^0(A,\Sym^r\Omega_A)
\end{eqnarray*}}
\end{prop}

\begin{proof} We follow the ideas of \cite{schn}. The symmetric  powers of the
exact sequence 
$0\to N_{X/A}^*\to \Omega_A\vert_X\to \Omega_X\to 0
$  
yield, for each nonnegative $r$, a long exact sequence
$$\begin{matrix}
 0\to  \wedge^cN^*_{X/A}\otimes \Sym^{r-c}\Omega_A 
\to \cdots\to N^*_{X/A}\otimes \Sym^{r-1}\Omega_A  \to 
\Sym^r\Omega_A\vert_X\to \Sym^r\Omega_X \to 0 
\end{matrix} 
$$
By Le Potier's vanishing theorem (\cite{lp};  \cite{laz}, Remark 7.3.6),
$H^q(X,\wedge^iN^*_{X/A} )$ vanishes for $n-c-q> c-i$ and $i>0$. Since
$\Omega_A$ is trivial, we get, by an elementary homological algebra argument
(\cite{schn}, Lemma, p. 176),
$$H^q(X,\Ker (\Sym^r\Omega_A\vert_X\to \Sym^r\Omega_X))=0\quad\hbox{for all
}q\le 2n-c
$$
The proposition now follows from the fact that 
the restriction $H^q(A,\cO_A)\to H^q(X,\cO_X)$, hence also the restriction
$H^q(A,\Sym^r\Omega_A)\to H^q(X,\Sym^r\Omega_A\vert_X)$, is bijective for
$q\le n-2c$ (\cite{som2}).\end{proof}

 Sommese proved (\cite{som}, Proposition (1.7)) that for any
$k$-ample vector bundle $\cE$ on a projective variety $X$ and  
any coherent sheaf $\cF$ on $X$, 
$$H^q(X,\Sym^r\cE\otimes \cF)=0$$
for all $q>k$ and $r\gg 0$. Theorem \ref{va} and Proposition \ref{p15} therefore imply the following.

 \begin{coro}\label{coro14}
Let $X$ be the  intersection of $c$ sufficiently ample\footnote{To be
more precise, we need condition (\ref{*})  to be satisfied.} general  hypersurfaces in
an abelian variety
$A$ of dimension
$n$. We have\footnote{Recall $h^q(A,\Sym^r\Omega_A
)=h^q(A,\cO_A) \dim\Sym^r \Omega_{A,0} =\binom{n}{q}\binom{n+r-1}{n-1}$.}
$$
h^q(X,\Sym^r\Omega_X )\ \begin{cases}=0&\hbox{ for $q>\max\{n-2c,0\}$ and }r\gg0\\
=h^q(A,\Sym^r\Omega_A )&\hbox{ for $q<n-2c$ and }r\ge0\\
\ge h^q(A,\Sym^r\Omega_A )&\hbox{ for $q=n-2c$ and }r\ge0
\end{cases}
$$
\end{coro}

\begin{rema}\upshape\label{ree} Assume $X$ as above is the complete intersection in $A$ of hypersurfaces with
classes $\ell_1,\dots,\ell_c$. Let $d=n-c$ be the dimension of $X$. 

 If $2c\ge n$ and $r\gg0$, the only nonzero cohomology group of $ \Sym^r\Omega_X$ is  $H^0(X,\Sym^r\Omega_X )$. The Segre class $s_d (\Omega^*_X)$ is given by (\cite{fult},  Example 3.2.12)
$$s_d (\Omega^*_X)=\bigl[ 
\prod_{i=1}^c(1+\ell_i\vert_X)\bigr]_d\\ =\sum_{1\le i_1<\dots<i_d\le c } 
\ell_{i_1}\cdots\ell_{i_d}\cdot\ell_1\cdots
\ell_c
$$ 
and (\cite{fult},  \S 3.1, where, however,
$\P(\cE)$ is the projective bundle of {\em lines} in the fibers of $\cE$)
\begin{eqnarray*} h^0(X,\Sym^r\Omega_X)\!\!\!&=&\!\!\!
\chi(\P(\Omega_X),\cO_{\P(\Omega_X)}(r))= \frac{r^{2d-1}}{(2d-1)!}s_d
(\Omega^*_X)+O(r^{2d-2})\\ h^0(A,\Sym^r\Omega_A)\!\!\!&= &\!\!\!
\frac{r^{n-1}}{(n-1)!}+O(r^{n-2})
\end{eqnarray*}
 
For $2c>n$,  the restriction $H^0(A,\Sym^r\Omega_A )\to
H^0(X,\Sym^r\Omega_X )$ is therefore {\em not} injective for $r\gg0$. 

For
$2c=n$,  this restriction   is injective by Proposition \ref{p15}, but {\em
not} surjective   for $r\gg0$, because $s_d
(\Omega^*_X)=(\ell_1\cdots\ell_d)^2\ge n!>1$.

If now $2c< n$, the only nonzero cohomology groups of $\Sym^r\Omega_X $  are, for $r\gg0$, 
\begin{eqnarray*}
H^0(X,\Sym^r\Omega_X)&\isom& H^0(A,\Sym^r\Omega_A)\\
&\vdots\\
H^{n-2c-1}(X,\Sym^r\Omega_X)&\isom& H^{n-2c-1}(A,\Sym^r\Omega_A)\\
H^{n-2c}(X,\Sym^r\Omega_X)&\supset& H^{n-2c}(A,\Sym^r\Omega_A)
\end{eqnarray*}
For $c>0$, direct calculations show that the latter inclusion is strict.
\end{rema}

\section{Subvarieties of the projective space}\label{s2}

We now study the positivity properties of the cotangent bundle of a
smooth subvariety of the projective space.

\subsection{Big twisted cotangent bundle}\label{ss21}

If $X$ is a smooth subvariety of $\P^n$ of dimension $d$, we let $\gamma_X:X\to
G(d,\P^n)$ be the Gauss map. We denote by 
$\cS$ the  universal subbundle and by $\cQ$  the universal quotient bundle on $ G(d,\P^n)$.  We
have 
$\gamma_X^*\cQ=N_{X/\P^n}(-1)$ and a commutative diagram
\begin{equation}\label{dia}
\begin{matrix}
&&0&&0\\
&&\downarrow&&\downarrow\\
&&N_{X/\P^n}^*(1)&=&N_{X/\P^n}^*(1)\\
&&\downarrow&&\downarrow\\
0&\lra&\Omega_{\P^n}(1)\vert_X&\lra&\cO_X^{n+1}&\lra&\cO_X(1)&\lra&0\\
&&\downarrow&&\downarrow&&\Vert\\
0&\lra& \Omega_X(1)&\lra&\gamma_X^*\cS^*&\lra&\cO_X(1)&\lra&0\\
&&\downarrow&&\downarrow\\
&&0&&0\\ 
\end{matrix}
\end{equation}

 \begin{prop}\label{rbig} 
 Let $X$ be a smooth subvariety of dimension $d$ of $\P^n$. 
\begin{itemize}
\item If 
$\gamma_X^*\cS^*$ is   big,  $2d\le n$.
\item If $2d\le n$ and $N_{X/P^n}(-1) $ is ample,\footnote{This is equivalent
to the following condition: for any hyperplane $H$ in $\P^n$, the set $\{x\in
X\mid {\bf T}_{X,x}\subset H\}$ is finite.
It holds for smooth nondegenerate complete intersections.}   $
\gamma_X^*\cS^*$ is nef and big.
\end{itemize}
\end{prop}

It was
proved  in
\cite{schn}, by a different method, that if $\Omega_X(1)$ is  big,   $2d\le n$ (see Remark \ref{remo}(1)). This can also be obtained via the (well-known) argument we use below.

\begin{proof} 
The   analog of the map (\ref{f1}) is
\begin{equation}\label{f2}
\SelectTips{cm}{}\xymatrix@M=8pt@R=20pt@C=25pt{
\P(\gamma_X^*\cS^*)\ar@/^1.5pc/^f <11mm,3mm>;[rr] \ar@<-0.3ex>@{^{(}->}[r]&
\P^n\times X\ar[d]^-{p_2}
\ar[r]^-{p_1}&\P^n\\
&X }
\end{equation}
whose image is the tangential variety 
$\bigcup_{x\in X}{\bf T}_{X,x}$ of
$X$, where ${\bf T}_{X,x} $ is
the embedded tangent space to $X$ at $x$, a  linear subspace of $\P^n$
of dimension
$d$.
 The vector bundle $\gamma_X^*\cS^* $ is nef, and it is big \iff\   this
variety has the expected dimension
$2d$. 
This proves the first point.

 The second point follows from Lemma~\ref{le3} applied to the middle vertical
exact sequence of  diagram (\ref{dia}) above.\end{proof}

We apply the same ideas to prove an analog of Theorem \ref{va}.

\begin{theo}\label{pr}
Let $X$ be a general complete intersection in $\P^n$ of multidegree
$(e_1,\dots,e_c)$. If $e_1\ge 2$ and   
$e_2,\dots,e_c$ are all   
$ \ge n+2$, the vector bundle  
$\gamma^*_X\cS^*$ is $\max(n-2c,0)$-ample.
\end{theo}

\begin{proof}We need to prove that the fibers of the map $f$ in
(\ref{f2}) have dimension at most $m=\max(n-2c,0)$. This means that for $H_i$
general in $|\cO_{\P^n}(e_i)|$ and for {\em   any} 
$t$ in
$\P^n$, the dimension of the set of points $x$
in
$ X$ such that
$t\in {\bf T}_{X,x}$ is at most $m$. Pick coordinates and write
$t=(t_0,\dots,t_n)$. If $s$  is an equation of a hypersurface
$H$,  we let $\partial_t H$ be the hypersurface
with equation $\partial_t s=\sum_{i=0}^n t_i\frac{\partial s}{\partial
x_i}
$. With this notation, we want
$$\dim(H_1\cap \partial_t H_1\cap\dots\cap
H_c  \cap  \partial_t H_c)\le m$$ 
As in the proof of Theorem \ref{va}, we proceed by induction on $c$,
assuming $c\le n/2$. When $c=1$, it is clear that $e_1\ge2$ is
sufficient.

Assume
$Y_t=H_1\cap
\partial_tH_1\cap\dots\cap H_{c-1} 
\cap \partial_tH_{c-1}$ has (pure) codimension $2c-2$ in $\P^n$, with
irreducible components $Y_{t,1},\dots,Y_{t,m} $. Set
$Y=(Y_{t,i})_{\rm red} $, with normalization  $\nu:\widehat Y\to Y$; it follows from Lemma \ref{le}   that  
$  \nu^*H$ is integral of codimension $1$ in $\widehat Y$ for
$H$ outside a closed subset of codimension $\ge d-1 $ in $|\cO_{\P^n}(d)|$.

Assume that this is the case. If $\codim_Y(Y\cap H\cap \partial_t
H)\le 1$, the section $\nu^*(\partial_ts) $ must vanish on $  \nu^*H$. Since the
restriction
$$H^0(\widehat  Y, \cO_{\widehat Y}(d-1))\to H^0(\nu^*H,
\cO_{ \nu^*H}(d-1))
$$
is injective, it must also vanish on $\widehat Y$. Since any $d$ distinct points
of
$Y$ impose independent conditions on elements of $|\cO_{\P^n}(d-1)|$ and
the map
$\partial_t: H^0(\P^n, \cO_{\P^n}(d))\to H^0(\P^n, \cO_{\P^n}(d-1))
$
is surjective, we have proved that
the set of hypersurfaces
$H$  in $|\cO_{\P^n}(d)|$ such that
$\codim_{Y_t}(Y_t\cap H\cap \partial_t H)\le 1$ has codimension $\ge d-1$
in
$|\cO_{\P^n}(d)|$. The theorem follows. 
\end{proof}

 \begin{coro} 
Let $X$ be a general complete intersection in $\P^n$ of multidegree
$(e_1,\dots,e_c)$. If $e_1\ge 2$ and   
$e_2,\dots,e_c$ are all   
$ \ge n+2$, and $c\ge n/2$, the vector bundle  
$\Omega_X(1)$ is big.
\end{coro}

\begin{remas}\upshape\label{remo} (1)  For {\em any} smooth subvariety $X$ of   $\P^n$ of codimension $c<n/2$, we have by \cite{schn}, Theorem 1.1,
$$H^0(\P(\Omega_X),\cO_{\P(\Omega_X)}(r))\isom
H^0(X,\Sym^r(\Omega_X(1)) )=0$$ for
$r> 0$. In particular, {\em  $\Omega_X(1)$ is not big} and  {\em not  
$(n-2c)$-ample,} as Theorem
\ref{va}  could suggest.  

(2)   When $X$ is a {\em surface,} (\ie, $c=n-2$) results of Bogomolov (\cite{bog3}, \cite{bog4}) give the much better result that 
$ \Omega_X(-\frac15K_X)$ is big.
\end{remas}

\begin{proofcor} 
The last row of diagram (\ref{dia}) yields, for all positive integers $r$,
an exact sequence
$$0\to \Sym^r\bigl(
\Omega_X(1)\bigr) \to \Sym^r(\gamma^*_X\cS^*)\to \Sym^{r-1}(\gamma^*_X\cS^*)\otimes\cO_X(1)\to 0
$$
It follows from Theorem \ref{pr} that for $r\gg0$, we have 
\begin{equation}\label{hq}
H^q\bigl(X,  \Sym^r\bigl(
\Omega_X(1)\bigr)\bigr)=0
\quad\hbox{ for } q> 1 
\end{equation}
On the other hand, if $d=n-c$, the coefficient of $\frac{r^{2d-1}}{(2d-1)!}$ in
the polynomial
$\chi\bigl(X,  \Sym^r\bigl(
\Omega_X(1)\bigr)\bigr)$ is (Remark \ref{ree})
\begin{eqnarray*} 
s_d (\Omega_X(1)^*)&=& 
 \bigl[\prod_{i=1}^c(1+(e_i-1)h)(1-h)\bigr]_d\\
&=&\hskip -4mm\sum_{1\le i_1<\dots<i_d\le c}\hskip -4mm  (e_{i_1}-1)\cdots
(e_{i_d}-1)-\hskip -4mm\sum_{1\le i_1<\dots<i_{d-1}\le c}\hskip -6mm (e_{i_1}-1)\cdots
(e_{i_{d-1}}-1)
\end{eqnarray*}
Since $c\ge n/2$, this is  positive, so that by (\ref{hq}), we have, for $r\gg0$,
$$h^0\bigl(X,  \Sym^r\bigl(
\Omega_X(1)\bigr)\bigr)\ge\chi\bigl(X,  \Sym^r\bigl(
\Omega_X(1)\bigr)\bigr) =\a
r^{2d-1} +O(r^{2d-2})
$$
for some $\a>0$. This shows that $\Omega_X(1)$ is big.
\end{proofcor}

\subsection{Ample cotangent bundle}
 
By analogy with Theorem \ref{va}, it is tempting to conjecture the
following  generalization of a question formulated by
Schneider in 
\cite{schn}, p.\ 180.

\begin{conj}\label{conj11} The
cotangent bundle of the intersection in $\P^n$ of at least $n/2$ general
hypersurfaces of sufficiently high degrees is  ample.\footnote{For $n\ge 4$, this conjecture was recently proved by D. Brotbek for complete intersections of $n-2$ general  hypersurfaces of degrees $\ge\frac{8n+2}{n-3}$.}
\end{conj} 


We give in the next subsection alternative versions of this conjecture.
Note that   there are smooth complete intersections of arbitrarily large
multidegree whose cotangent bundle is not even nef. Here are two examples.

Let $X$ be a smooth subvariety of $\P^n$ that contains a line $\ell$.  The restriction $\Omega_X(2)\vert_\ell$ has the trivial quotient
$\Omega_\ell(2)$, hence is not ample. It follows that $\Omega_X(2)$ {\em is not ample.}
Note that there are complete intersections $X\subset \P^n$ of any dimension $\ge2$ and arbitrarily large degrees that contain a line.

Consider also 
the morphism $f:\P^n\to\P^n$ given by
$f(x_0,\dots,x_n)=(x_0^2,\dots,x_n^2)$. If $Y\subset\P^n$ is a general
complete intersection of multidegree
$(d_1,\dots,d_c)$, its inverse image 
$X=f^{-1}(Y)$ is a smooth complete intersection of multidegree
$(2d_1,\dots,2d_c)$. One checks that there is an exact sequence
$$0\lra (f^*\gamma_Y^*\cS^*)(-1)\lra \gamma_X^*\cS^*\lra  
\bigoplus_{i=0}^n\cO_{ X\cap H_i} \lra 0
$$
that implies that $\gamma_X^*\cS^*$, and {\it a fortiori} $\Omega_X(1) $, are not ample when $n-c\ge2$.


\subsection{Cohomology of symmetric tensors}   

Let $X$ be a general complete intersection in $\P^n$ of multidegree
$(e_1,\dots,e_c)$, with $e_1\ge 2$ and 
$e_2,\dots,e_c \ge n+2$. With the notation of  \S\ref{ss21}, it follows
from Theorem \ref{pr} and the result of Sommese used in \S\ref{ss13} that
for any coherent sheaf $\cF$ on $X$, 
\begin{equation}\label{som} H^q(X, \Sym^r(\gamma^*_X\cS^*)\otimes\cF)=0
\end{equation} for all $q>q_0= \max\{n-2c,0\}$ and $r\gg0$.

For ordinary ampleness, there is a converse to  Sommese's result
(\cite{laz}, Theorem 6.1.10).

\begin{prop} 
Let $X$ be a projective variety, let $\cE$ be a vector bundle on $X$, and let
$L$ be an ample line bundle on $X$. The following properties are equivalent:
\begin{itemize}
\item[\rm (i)] $\cE$ is  ample;
\item[\rm (ii)] for any integer $m$, we have 
$H^q(X,\Sym^r\cE\otimes L^m)=0$
for all $q>0$ and $r\gg 0$.
\end{itemize}
\end{prop}

\begin{proof}  Let $\cF$ be an arbitrary coherent sheaf on
$X$. It has a possibly nonterminating resolution
$$\cdots\to\cE_2\to\cE_1\to\cE_0\to\cF\to0$$
by locally free sheaves   that are direct sums of powers of $L$.
Therefore, we have $H^q(X,\Sym^r\cE\otimes\cE_j)=0$ for all
$j\in\{0,\dots,\dim(X)\}$, all
$q>0$ and $r\gg 0$, and this implies $H^q(X,\Sym^r\cE\otimes\cF)=0$ for   all
$q>0$ and $r\gg 0$. This is the usual cohomological criterion for the  
ampleness of $\cE$.\end{proof}

  Conjecture \ref{conj11} is therefore equivalent to the
following cohomological statement.

\begin{conj}\label{conj16}
Let $X$ be the intersection   in $\P^n$ of a least
$n/2$ general hypersurfaces of sufficiently high degrees. For any integer $m$, we have
\begin{equation}\label{van1}
H^q(X,(\Sym^r\Omega_X)(m))=0
\end{equation}
for all $q>0$ and $r\gg 0$.
\end{conj}

If $X$ is a smooth projective variety of dimension $d$
 with
$\omega_X$ nef and big, and $r>d$, we have (\cite{bog3} for $d=2$; \cite{dem}, Theorem
14.1 in general)
$$H^d(X,\Sym^r\Omega_X)\isom
H^0(X,\Sym^rT_X\otimes\omega_X)^*=H^0(X,\Gamma^{(r-1,-1,\dots,-1)}T_X
)^*=0$$
Assume $\omega_X$ is ample. Given any line bundle $L$ on $X$, there is a
positive integer
$r_0$ such that $\omega_X^{r_0-1}\otimes L$ is ample. For $r\ge r_0d$, we
have
\begin{eqnarray*}
H^d(X,\Sym^r\Omega_X\otimes L)&\isom&
H^0(X,\Sym^rT_X\otimes\omega_X\otimes
L^*)^*\\
&\isom&
H^0\bigl(X,\Sym^rT_X\otimes\omega_X^{r_0}\otimes(\omega_X^{r_0-1} \otimes
L)^*\bigr)^*\\
&=&H^0\bigl(X,\Gamma^{(r-r_0,-r_0,\dots,-r_0)}T_X\otimes(\omega_X^{r_0-1} \otimes
L)^* \bigr)^*=0
\end{eqnarray*}
by the theorem of Demailly mentioned above.
This leads to think that the following stronger form   of
Conjecture
\ref{conj16} might be true (compare with Corollary \ref{coro14}).

 \begin{conj}\label{conj17}
Let $X$ be the  intersection   in $\P^n$ of $c$ general
  hypersurfaces of sufficiently high degrees and let $m$  be an integer. For
$r\gg 0$, we have
\begin{equation}\label{van2} 
H^q(X,(\Sym^r\Omega_X)(m))=0
\end{equation}
except for $q=\max\{n-2c,0\}$.
\end{conj}

\begin{remas}\upshape (1)  For {\em any} smooth subvariety $X$ of   $\P^n$ of codimension $c$, the vanishing (\ref{van2}) 
holds for   $q<n-2c$ and $r\ge m+2$ by \cite{schn}, Theorem 1.1, and for
$q=n-c$ by Demailly's theorem. In particular, Conjecture
\ref{conj17} holds for  
 $c\le1$.  

(2) Recall (Remark \ref{ree}) that if $X$ is a smooth projective variety
of dimension $d$, the coefficient of $\frac{r^{2d-1}}{(2d-1)!}$ in the
polynomial
$\chi(\P(\Omega_X),\cO_{\P(\Omega_X)}(r))$ is
$s_d(\Omega^*_X)$. When $X$ is a smooth complete
intersection   of multidegree
$(e_1,\dots,e_c)$ in $\P^n$, we have (\cite{fult},
Example 3.2.12)
\begin{eqnarray*}
s_d (\Omega^*_X)&=&\left[ 
\frac{\prod_{i=1}^c(1+e_ih)}{(1+h)^{n+1}}\right]_d\\
&=&\sum_{1\le i_1<\dots<i_r\le c\atop 0\le r\le d} (-1)^{d-r}e_{i_1}\cdots
e_{i_r}\binom{n+d-r}{n}
\end{eqnarray*}
This is a polynomial in $ e_1,\dots,e_c$ whose leading term is, 
for $ 2c\le n$,
$$(-1)^{n-2c}e_1\cdots e_c\binom{2n-2c}{n}$$
and, for $2c\ge n$,
$$ \sum_{1\le i_1<\dots<i_d\le c}e_{i_1}\cdots e_{i_d}$$
When
$e_1,\dots,e_c
$ are big enough, its sign is therefore $(-1)^{\max\{n-2c,0\}}$. This is
compatible with Conjecture \ref{conj17}.
\end{remas}

\section{Bogomolov's construction of varieties with ample cotangent
bundle}\label{s4}

We present here a construction due to Bogomolov that produces varieties with
ample cotangent bundle as linear sections of products of varieties with big
cotangent bundle. Bogomolov's  construction
appears in \cite{wo}  in a differ\-ential-geometric setting. Everything in
this section is due to Bogomolov.\footnote{I am grateful to  Bogomolov for 
allowing me to reproduce his construction.}

  \subsection{The construction}

 \begin{prop}[Bogomolov]
 Let
$X_1,\ldots,X_m$ be smooth projective varieties   with 
big cotangent bundle, all of dimension at least $d>0$. Let $V$ be a general
linear section of $X_1\times\dots\times X_m$.
If $\dim(V)\le \frac{d(m+1)+1}{2(d+1)}$, the cotangent bundle of $V$ is
ample.
\end{prop}

A variant of this construction appears in \cite{laz}, 6.3.34. The bound
that Lazarsfeld's argument yields, $\dim(V)\le\frac{d(m+1)}{2d+1}$, is slightly
better, but one needs to take hyperplane sections of sufficiently high
degree.

\begin{proof}
Since $\Omega_{X_i}$ is big, there exist a proper closed subset $
B_i$ of
$\P(\Omega_{X_i})$ and an integer $q$ such that for each $i$,  the sections of
$\cO_{\P(\Omega_{X_i})}(q)$, i.e., the sections of $\Sym^q\Omega_{X_i}$, define
an {\em injective} morphism
$$f_i:\P(\Omega_{X_i})\moins B_i\lra \P^{n_i} 
$$
 
\begin{lemm} 
 Let
$X $ be a smooth subvariety of a projective space and let $B$ be a  
subvariety  of $\P(\Omega_X)$. A general linear
 section $V$ of
$X$ of dimension at most $ \frac12\codim(B)$ satisfies
$$\P(\Omega_V)  \cap  B  = \vide $$
\end{lemm}

\begin{proof} Let $\P^n$ be the ambiant projective space. Consider the variety
$$\{((t,x),\Lambda)\in B\times G(n-c,\P^n)\mid x\in X\cap \Lambda,\ t\in 
T_{X,x}\cap T_{\Lambda,x}   )
$$ The fibers of its projection  to $B $ have codimension $2c$,
hence it does not dominate $G(n-c,\P^n) $ as soon as $2c>\dim(B)$. This is
equivalent to $2(\dim(X)-\dim(V))-1\ge 2\dim(X)-1-\codim(B)$ and   the
 lemma is proved.
\end{proof}

 Let $B'_i$ be the (conical) inverse image
of
$B_i$ in the total space  of the tangent bundle of $X_i$. Let $V$ be a
general linear section of
$X_1\times\dots\times X_m$ and set $a=m+1-2\dim(V)$.

 If 
$t=(t_1,\dots,t_m)$, with
$t_i\in T_{X_i,x_i}$, is a nonzero tangent vector to
$V$, the lemma  implies that there are at least
$a$ values of the index $i$ for which $t_i\notin B'_i$. If, say, $t_1$ is not in
$B'_1$, there exists a section of $\Sym^q\Omega_{X_1}$ that does not
vanish at $t_1$. This section induces, via the projection $ V\to X_1$, a  
section of
$\Sym^q\Omega_V$ that does not vanish at $t$. It follows that
$\cO_{\P(\Omega_V)}(q)  $ is base-point-free and its sections define
a morphism
 $f:\P(\Omega_V) \lra \P^n
$.

We need to show that $f$ is {\em finite}. Assume to the contrary  that a
curve $C$ in $
\P(\Omega_V)$ through
$t$ is contracted. Since the restriction of  the projection $\pi:\P( \Omega_V
)\to V$ to any fiber of $f$ 
 is injective, and since $f_i$ is injective, the
argument above proves that the curve
$\pi(C)$ is contracted by each projection
$p_i:V\to X_i$ such that $t_i\notin B'_i$. 

The following lemma leads to a contradiction when $2\dim(V)\le  ad+1$. This
proves the proposition.
\end{proof}

\begin{lemm} 
 Let $V$ be a general linear section of a product $X\times Y$ in a projective
space. If $2\dim(V)\le \dim(X)+1$, the projection $V\to X$ is finite.
\end{lemm}

\begin{proof} Let $\P^n$ be the ambiant projective space. Consider the
variety that is the closure of
$$\{(x,y,y',\Lambda)\in X\times Y\times Y\times
G(n-c,\P^n)\mid  y\ne y',\ (x,y)\in 
\Lambda,\ (x,y')\in \Lambda  )\}
$$ The fibers of its projection  to $X\times Y\times Y$ have codimension
$2c$, hence the general fiber of its projection to $G(n-c,\P^n) $, which is the
closure of
$$\overline{\{((x,y),(x,y'))\in V\times V\mid \\ y\ne y',\
(x,y)\in  V,\ (x,y')\in V  )\}}$$ has dimension at most $1$ as soon as
$2c\ge \dim(X\times Y\times Y)-1$. This implies that the projection $V\to
X$ is finite and is equivalent to
$2\dim(V)\le \dim(X)+1$:   the
 lemma is proved.
\end{proof}

 \subsection{Fundamental groups}

Using this construction, Bogomolov exhibits smooth projective  varieties with
ample cotangent bundle that are simply connected. More generally, his ideas give the following result.

 \begin{prop}
Given any smooth projective variety $X$, there exists a\break smooth projective
surface with ample cotangent bundle and same fundamental
group as $X$.
\end{prop}

\begin{proof}
By the Lefschetz hyperplane theorem, a sufficiently ample $3$-dimen\-sional
linear section $Y$ of $X\times\P^3$ has same fundamental group as $X$ and
$K_Y$ ample. A smooth hyperplane section $S$ of $Y$ with  class $ah$ satisfies
$$c_1^2(S)-c_2(S)=  a^2h^2\cdot c_1(Y)+ah\cdot (c_1^2(Y)-c_2(Y)) 
$$  This is positive for $a\gg0$. By a famous trick of Bogomolov (\cite{bog4}),
$S$ is a smooth surface of general type with big cotangent bundle and same
fundamental group as
$X$. Starting from a simply connected $X_0$, we similarly obtain a simply
connected surface $S_0$ with big cotangent bundle. Taking in Bogomolov's
construction
$X_1=\dots=X_m=S_0$, we produce a smooth simply connected projective surface
$S_1$ with ample cotangent bundle. 

Taking in Bogomolov's construction
$X_1=S$ and
$X_2=\dots=X_m=S_1$, we produce a smooth projective surface with ample
cotangent bundle and same fundamental group as
$X$.
\end{proof}

\subsection{Geography}

At the end of the last century, it was traditional to look at the so-called geography of 
Chern numbers of surfaces of general type with a given property. For the
sake of the old days, let us look at the numbers $c_1^2$ and $c_2$ for a
surface  with ample cotangent bundle. From \cite{miy} and \cite{fula}, we
get
$$3c_2 \ge c_1^2>c_2
$$ All surfaces of general type with  $3c_2 = c_1^2 $  have ample
cotangent bundle (\cite{mi}, Corollary, p. 294). Examples of such surfaces
were constructed in \cite{hir}. 

Let 
$T$ be a threefold with ample cotangent bundle. A smooth hyperplane section
$S$ of $T$ with  class $ah$ has ample cotangent bundle. It satisfies
$$\frac{c_1^2(S)}{c_2(S)} = \frac{a^2h^3+2ah^2\cdot c_1(T)+ah\cdot
c_1^2(T)}{a^2h^3+ah^2\cdot c_1(T)+ah\cdot c_2(T)}  
$$  and this ratio is as close to $1$ as we want for $a\gg0$. It is likely
that using Bogomolov's construction and spending some time on tedious
calculations, one could produce smooth projective surfaces with ample
cotangent bundle such that the ratio $c_1^2(S)/c_2(S)$ is any given rational
number in $(1,3]$ (at least, these ratios should fill out a dense open subset of this interval).

\end{document}